\documentclass{commat}

\title[Existence of normal $p$-complement of finite groups with restrictions]{On existence of normal $p$-complement of finite groups with restrictions on the conjugacy class sizes}

\author{Ilya B. Gorshkov}

\affiliation{
\address{Sobolev Institute of Mathematics SB RAS
Novosibirsk, Russia
Siberian Federal University,
Krasnoyarks, Russia}

\email{ilygor8@gmail.com}
}

\thanks{
    The work is supported by Russian Science Foundation (project 18-71-10007).
}

\abstract{The greatest power of a prime $p$ dividing the natural number $n$ will be denoted by $n_p$. Let $Ind_G(g)=|G:C_G(g)|$. Suppose that $G$ is a finite group and $p$ is a prime. We prove that if there exists an integer $\alpha>0$ such that $Ind_G(a)_p\in \{1,p^{\alpha}\}$ for every $a$ of $G$ and a $p$-element $x\in G$ such that $Ind_G(x)_p>1$, then $G$ includes a normal $p$-complement.}



\firstpage{93}

\VOLUME{30}

\DOI{https://doi.org/10.46298/cm.9294}

\begin{paper}

\section{Introduction}
In this paper, all groups are finite. Denote the set of prime divisors of positive integer $n$ by $\pi(n)$, and by the set $\pi(|G|)$ for a group $G$ by $\pi(G)$. For a set of primes $\pi$ and a positive integer $n$ we will denote $n_{\pi}=\prod_{p\in\pi}|n|_p$. Let $G$ be a group and take $a\in G$. With $a^G$ standing for the conjugacy class in $G$ containing $a$, put $N(G)=\{|x^G|, x \in G\} \setminus\{1\}$. Denote by $|G||_p$ the number $p^n$ such that $N(G)$ contains a multiple of $p^n$ and avoids multiples of $p^{n+1}$. For $\pi\subseteq\pi(G)$ put $|G||_{\pi}=\prod_{p\in \pi}|G||_p$. For brevity, $|G||$ is meaning $|G||_{\pi (G)}$. Observe that $|G||_p$ divides $|G|_p$ for each $p\in\pi(G)$. However, $|G||_p $ can be less than $|G|_p$. 
Take a set of primes $\pi$, denote $\Theta_{\pi}=\{\tau\subseteq\pi\ |\ \tau\neq \varnothing, |\tau|\geq|\pi|-1\}$.
\begin{definition}
Let $\pi$ be a set of primes. We say that a group $G$ satisfies the condition $\pi^*$ or $G$ is a $\pi^*$-group and write $G\in\pi^*$ if for every $a\in N(G)$ there exists $\tau_{a}\in \Theta_{\pi}$ such that $a_{\pi}=|G||_{\tau_{a}}$.
\end{definition}
 Ishikawa \cite{Ishikawa} proved that a group $G$ with $N(G)=\{p^{\alpha}\}$ is nilpotent class at most $3$. Casolo, Dolfi and Jabara \cite{Casolo} described the set of $\{p\}^*$-groups. In particular, they proved that any group of $\{p\}^*$ is solvable and includes a normal $p$-compliment.
Camina \cite{Camina} proved that a group $G$ with $\{p,q\}^*$-property is nilpotent if $N(G)=\{p^n,q^m,p^nq^m\}$. Beltram and Filipe \cite{Beltran} extended Camina's theorem in the following way. Let $G$ be a group whose set of conjugacy class sizes is $\{1,n,m,nm \}$, where $n$ and $m$ are coprime positive integers; then $G$ is nilpotent and the integers $n$ and $m$ are prime-power numbers; in particular $G\in\{n,m\}^*$. The author of \cite{GorshkovC} investigate $\{p,q\}^*$-groups with trivial center.

In the present paper we will investigate some generalisations of the property $\{p\}^*$.

\begin{definition}
We say that a group $G$ satisfies the condition $R(p)$ or $G$ is a $R(p)$-group and write $G\in R(p)$ if there exists an integer $\alpha>0$ such that $a_p\in\{1,p^{\alpha}\}$ for each $a\in N(G)$.
\end{definition}
Note that, if $G\in \pi^*$, then $G\in R(p)$ for each $p\in \pi$. The set of $R(p)$-group disjoins on two subsets $R(p)^*$ and $R(p)^{**}$:

\begin{enumerate}
\item[$(i)$]{ $G\in R(p)^*$ if $G$ contains a $p$-element $h$ such that $Ind_G(h)_p>1$;}

\item[$(ii)$]{ $G\in R(p)^{**}$ if $Ind_G(h)_p=1$ for each $p$-element $h\in G$.}
\end{enumerate}

We prove the following theorem.

\begin{theorem}
If $G\in R(p)^*$, then $G$ has a normal $p$-complement.
\end{theorem}
It follows from the theorem that the center of a $R(p)^*$-group is not trivial.

\begin{corollary}
If $G\in R(p)^*$ and $P\in Syl_p(G)$, then $Z(P)\leq Z(G)$.
\end{corollary}

Vasil'ev \cite{Vasilev} proved that if $G$ is a $R(p)$-group with trivial center and $|G||_p=p$, then Sylow $p$-subgroups of $G$ are abelian. This assertion is true in the general case.

\begin{corollary}
If $G\in R(p)$ and $Z(G)=1$, then Sylow $p$-subgroups of $G$ are abelian.
\end{corollary}

\section{Preliminary results}

\begin{lemma}[{\rm \cite[Lemma 1.4]{GorshkovA}}]\label{factorKh}
For a finite group $G$, take $K\unlhd G$ and put $\overline{G}= G/K$. Take $x\in G$ and $\overline{x}=xK\in G/K$.
The following claims hold:
\begin{enumerate}
\item[(i)] $|x^K|$ and $|\overline{x}^{\overline{G}}|$ divide $|x^G|$.

\item[(ii)] For neighboring members $L$ and $M$ of a composition series of $G$, with $L<M$, take $x\in M$  and
the image $\widetilde{x}=xL$ of $x$. Then $|\widetilde{x}^S|$ divides $|x^G|$, where $S=M/L$.

\item[(iii)] If $y\in G$ with $xy=yx$ and $(|x|,|y|)=1$, then $C_G(xy)=C_G(x)\cap C_G(y)$.

\item[(iv)] If $(|x|, |K|) = 1$, then $C_{\overline{G}}(\overline{x}) = C_G(x)K/K$.

\item[(v)] $\overline{C_G(x)}\leq C_{\overline{G}}(\overline{x})$.
\end{enumerate}
\end{lemma}

\begin{lemma}[{\rm \cite[Lemma 2.1]{Casolo}}]\label{Dol}
Let $x,y$ be elements of a group $G$ and assume at least one of the following conditions:
\begin{enumerate}
    \item[(i)] $x$ and $y$ commute and have coprime orders;
    \item[(ii)] $x\in N, y\in M$ with $N,M \unlhd G$ and $N\cap M=1$.
\end{enumerate}
Then $C_G(xy)=C_G(x)\cap C_G(y)$.
\end{lemma}

\begin{lemma}[{\rm \cite[Lemma 2.7]{Casolo}}]\label{Dol2}
Let $A$ be a group acting via automorphisms on a group $G$ and $N$ be a normal $A$-invariant subgroup of $G$.  If $(|A|,|N|)=1$, then:
 \begin{itemize}
     \item [(i)] $C_{G/N}(A)=C_G(A)N/N$;
     \item [(ii)] $|C_G(A)|=|C_N(A)||C_{G/N}(A)|$.
 \end{itemize}
 \end{lemma}
 \begin{lemma}[{\rm \cite[Lemma 4]{GorshkovA}}]\label{hh}
Take $g\in G$. If each conjugacy class of $G$ contains an element $h$ such that $g\in C_G(h)$, then $g\in Z(G)$.
\end{lemma}

\begin{lemma}\label{norm}
If $G\in R(p)$ and $N\unlhd G$ such that $|N|_p=|G|_p$, then $N\in R(p)$ or $|N||_p=1$.
\end{lemma}
\begin{proof}
Since $N$ is a normal subgroup and $N$ includes every Sylow $p$-subgroup of $G$, we have $N$ includes every Sylow $p$-subgroup of $C_G(x)$ for any $x\in G$. Therefore, $Ind_G(x)_p=Ind_N(x)_p$ and the lemma is proved.
\end{proof}

\begin{lemma}\label{pshtrih}
If $G\in R(p)$, $N\unlhd G$ is a $p'$-group, then $G/N\in R(p)$ or $|G/N||_p=1$.
\end{lemma}
\begin{proof}
Let $\overline{\textcolor{white}{a} }: G\rightarrow G/N$ be a natural homomorphism. We have $Ind_G(h)$ is a multiple of $Ind_{\overline{G}}(\overline{h})$ for any $h\in G$. Therefore, $|G||_p\geq |\overline{G}||_p$. Assume that there exists $\overline{x}\in \overline{G}$ such that $1<Ind_{\overline{G}}(\overline{x})_p<|G||_p$. Let $H$ be a Sylow  $p$-subgroup of $C_{\overline{G}}(\overline{x})$. Therefore, $|G|_p>|H|>|G|_p/|G||_p$. Put $T<G$ is a $p$-group such that $\overline{T}=H$. From Lemma \ref{Dol2} follows that $C_G(T)$ contains $y$ such that $yN=\overline{x}$. Since $C_G(y)\geq T$, we obtain $Ind_G(y)_p\leq|G|_p/|T|<|G||_p$. Therefore, $Ind_G(y)_p=1$ and consequently $Ind_{\overline{G}}(yN)_p=1=Ind_{\overline{G}}(x)_p$; a contradiction.
\end{proof}
The prime graph $GK(G)$ of a finite group $G$ is defined as follows. The vertex set of $GK(G)$ is
the set $\pi(G)$. Two distinct primes $p,q\in\pi(G)$
considered as vertices of the graph are adjacent by the edge if and only if there is an element of
order $pq$ in $G$. Denote by $s(G)$ the number of connected components of $GK(G)$ and by $\pi_i(G), i=1,...,s(G)$, its $i$-th connected
component. If $G$ has even order, then put $2\in \pi_1(G)$.

\begin{lemma}\cite[Theorem A]{Williams}\label{GK}
If a finite group $G$ has disconnected prime graph, then one of the
following conditions holds:
\begin{enumerate}
\item[(a)]{$s(G)=2$ and $G$ is a Frobenius or 2-Frobenius group;}
\item[(b)]{there is a nonabelian simple group $S$ such that $S\leq G = G/F(G)\leq Aut(S)$, where $F(G)$ is the
maximal normal nilpotent subgroup of $G$; moreover, $F(G)$ and $G/S$ are $\pi_1(G)$-subgroups, $s(S)\geq s(G)$,
and for every $i$ with $2\leq i\leq s(G)$ there is $j$ with $2\leq j\leq s(S)$ such that $\pi_i(G)=\pi_j (S)$.}
\end{enumerate}
\end{lemma}
\begin{lemma}[{\rm \cite[Lemma 5.3.4]{Gorenstein}}]\label{Thom}
Let $A\times B$ be a group of automorphisms of the $p$-group $P$ with $A$ a $p'$-group and $B$ a $p$-group. If $A$ acts trivially on $C_P(B)$, then $A=1$.
\end{lemma}

\begin{lemma}[{\rm \cite[Lemma 5.2.3]{Gorenstein}}]\label{hz}
Let $A$ be a $p'$-group of automorphisms of the abelian group $P$. Then we have $P=C_P(A)\times [P,A]$
\end{lemma}

\begin{lemma}\cite[Lemma 1]{Camina}\label{Camina}
If, for some prime $p$, every $p'$-element of a group $G$ has index prime to $p$, then the Sylow $p$-subgroup of $G$ is a direct factor of $G$.
\end{lemma}
\begin{lemma}{\rm\cite[Lemma 3.6]{VVasilev}}\label{l:hallhigman} For distinct primes $s$ and $r$, consider
a semidirect product $H$ of a normal $r$-subgroup $T$ and a cyclic
subgroup $C=\langle g\rangle$ of order~$s$ with $[T,g]\neq1$.
Suppose that $H$ acts faithfully on a vector space $V$ of positive characteristic
$t$ not equal to~$r$. If the minimal polynomial of $g$ on $V$ does not equal $x^s-1$, then
\begin{enumerate}
\item[(i)] $C_T(g)\neq1$;
\item[(ii)] $T$ is nonabelian;
\item[(iii)] $r=2$ and $s$ is a Fermat prime.
\end{enumerate}
\end{lemma}

\begin{lemma}\cite[Lemma 11]{GorshkovB}\label{Simple}
If $S\leq A\leq Aut(S)$, where $S$ is a nonabelian simple group, then $|A|=|A||$.
\end{lemma}

\begin{lemma}\cite[Theorem B]{Navarro}\label{Navar}
Let $G$ be a finite group and $p$ a prime. Suppose that for every $p$-element $x$ the number $|x^G|$ is a $p'$-number.
Then,
\[
O^{p'}(G/O_{p'}(G)) = S_1 \times \dotsm \times S_r \times H,
\]
where $H$ has an abelian Sylow $p$-subgroup, $r\geq 0$ , and $S_i$ is a nonabelian simple group
with either
\begin{enumerate}
\item[(i)] $p=3$ and $S_i\simeq Ru$ or $J_4$ or $S_i\simeq\ ^2F_4(q)'$, $9\nmid (q+1)$; or

\item[(ii)] $p=5$ and $S_i\simeq Th$ for all $i$.
\end{enumerate}
\end{lemma}

\section{Proof}

Let $G$ be a counterexample for assertion of the theorem of minimal order.

\begin{lemma}\label{ops}
$O_{p'}(G)=1$
\end{lemma}
\begin{proof}
From Lemma \ref{pshtrih} it follows that $G/O_{p'}(G)\in R(p)$ or $|G/O_{p'}(G)||_p=1$. We can think that $G/O_{p'}(G)$ does not include a normal $p$-complement, else $G$ contains a normal $p$-complement. Therefore, $G/O_{p'}$ a counterexample for assertion of theorem; a contradiction with minimality $G$. If $|G/O_{p'}(G)||_p=1$, then Lemma \ref{Camina} implies that $G/O_{p'}(G)$ is a $p$-group. Therefore, $O_{p'}(G)$ is a normal $p$-compliment of $G$; a contradiction.
\end{proof}

\begin{lemma}\label{ops2}
Orders of minimal normal subgroups of $G$ are multiples of $p$.
\end{lemma}
\begin{proof}
It follows from Lemma \ref{ops}.
\end{proof}

\begin{lemma}\label{cokp}
Each minimal normal subgroup of $G$ is a $p$-group.
\end{lemma}
\begin{proof}
Let $H$ be the socle of $G$. Then $H$ has expression in form $S\times X$, where $S=S_1\times S_2\times...\times S_n$, for nonabelian simple groups $S_1,...,S_n$ and a $p$-group $X$.
It follows from Lemma \ref{ops2} that $p$ divides the order of $S_i$ for all $1\leq i\leq n$. Assume that $G$ contains a $p$-element $x$ such that $S_1^x\neq S_1$. Let $D=\langle aa^xa^{x^2}...a^{|x|-1}| a\in S_1\rangle$. We have $D=C_{S_1\times S^x_1...}(x)$ and $D\simeq S_1$. Since $|S_1|$ is a multiple of $p$, we see that $Ind_G(x)_p>1$. By Lemma \ref{Simple}, we obtain $D$ contains a $p'$-element $y$ such that $Ind_D(y)_p=|D|_p$. Thus, $Ind_{G}(xy)_p>Ind_{G}(x)_p=|G||_p$; a contradiction.
It follows that $S_i^y=S_i$, for any $1\leq i\leq n$ and for any $p$-element $y$. Take $h_1\in S_1$ and $h_2\in S_2$. We have a $p$-element $y\in C_G(h_1h_2)$ iff $y\in C_G(h_1)\cap C_G(h_2)$. Assume that $n>1$. Since $|S_i|=|S_i||$, we see that $S_i$ contains an element $h_i$ such that $|h_i^{S_{i}}|_p=|S_i|_p$, in particular $|h_i^G|_p>1$. Let $A$ be a Sylow $p$-subgroup of $C_G(h_1h_2)$ and $B$ be a Sylow $p$-subgroup of $C_G(h_1)$. Then $A<C_G(h_1)\cap C_G(h_2)$ and $|G|_p>|B|>|A|$; a contradiction. Therefore, $n=1$.

From Lemma \ref{Simple}, it follows that $S$ contains an element $h$ such that a $p$-element $y\in C_G(h)$ iff $y\in C_G(S)$. If $C_G(S)$ contains $y$ such that $|y^{C_G(S)}|_p>1$, then $|(yh)^G|_p>|h^G|_p$; a contradiction. From Lemma \ref{Camina} it follows that $C_G(S)=O_p(G)$. Moreover $O=O_p(G)$ is abelian. We have $|h^G|_p=|G|_p/|O|_p$. Take $a\in S$ such that $|a^S|_p<|S|_p$. Hence $|a^G|_p<|h^G|_p$. This implies that $|a^S|_p=1$. From Lemma \ref{Navar} it follows that a Sylow $p$-subgroup of $S$ is abelian or $S$ is isomorphic to one of groups $J_4, Ru,\ ^2F_4(q)', Th$. Also it signifies that $S\in R(p)$.

Assume that there exists a $p$-element $x\in G\setminus H$ such that $Ind_G(x)_p>1$ and $x$ acts on $S$ as an outer automorphism. From Lemma \ref{Camina} and the equation $Ind_{C_G(x)}(y)_p=1$ for each $p'$-element $y\in C_G(x)$ it follows that $C_G(x)=L\times T$, where $L$ is a Sylow $p$-subgroup of $C_G(x)$. Therefore, $C_{S}(x)=\widetilde{P}\times\widetilde{L}$, where $\widetilde{P}$ is a Sylow $p$-subgroup of $C_{S}(x)$.

Assume that $S\simeq Alt_n$, where $n\geq5$. Since $\pi(Out(S))=\{2\}$, we obtain $p=2$. If $n\neq 6$, then $C_{S_1}(x)\simeq Alt_{n-2}$; a contradiction. If $n=6$, then $C_{S}(x)\simeq Alt_{4}$ or $C_{S_1}(x)\simeq Sym_{3}$; a contradiction.

From \cite{Conway} it follows that $S$ is not isomorphic to a sporadic simple group or the Tits group.

Therefore, $S$ is a group of Lie type.
Assume that a Sylow $p$-subgroup of $S$ is nonabelian. Since $S$ is not isomorphic to one of a sporadic groups, it follows that $S\simeq \ ^2F_4(q)', 9\nmid q+1$ and $p=3$. Therefore, $x$ acts on $S$ as a field automorphism. Thus, $q=2^{3(2m+1)}$; this contradicts with that $9\nmid q+1$. Thus, Sylow $p$-subgroups of $S$ are abelian.

Assume that $p=2$. From description of simple groups with abelian Sylow $2$-subgroup \cite{Walter} it follows that $S$ is isomorphic to one of a groups $L_2(q)$ where $q=2^f$ or $q\equiv 3,5(mod\ 8)$, $J_1$ or $^2G_2(q)$  where $q=3^{2m+1}$ and $m\geq1$. Put $P\in Syl_2(S)$. From \cite[Theorems 1 and 7]{Kondrat'ev} it follows that $C_S(P)=Z(P)$ or $S$ is isomorphic to $L_2(q)$ where $q$ odd. Let $S$ be isomorphic to $L_2(q)$ for some odd $q$. If $x$ acts on $S$ as a field automorphism, then $C_S(x)\simeq L_2(d)$, where $d$ divides $q$, in particular $C_S(x)$ is not a direct product of a Sylow $p$-subgroup and a $p$-complement; a contradiction. Assume that $x$ acts on $S$ as diagonal automorphism or a diagonal-field automorphism. Therefore, $S$ contains a $2$-element $z$ such that $C_G(z)\cap x^G=\varnothing$. Consequently $1<|z^G|_p<|h^G|_p$; a contradiction. Hence, $S$ does not isomorphic $L_2(q)$ for odd $q$.
It follows that $C_S(P)=Z(P)$. Since $|S|=|S||$ and $S\in R(p)$, we see that $\{p\}$ is a connected component of $GK(S)$. The group of outer automorphisms of $J_1$ is trivial, therefore $S$ does not isomorphic $J_1$.
If $S\simeq ^2G_2(q)$, then from \cite{VasilevVdovin} it follows that $2$ is not a connected components of $GK(S)$. If $S\simeq L_2(q)$ for even $q$, then $Out(S)$ is isomorphic to a group of field automorphisms. By analogy as before we can assume that $C_S(x)$ is not a direct product of a Sylow $p$-subgroup and a $p$-complement; a contradiction.
Thus, $p>2$.

From description of finite simple groups with an abelian Sylow $p$-subgroup \cite{Shen} it follows that $p$ does not divide the orders of graph and diagonal automorphism groups. Lemma \ref{norm} and fact that subgroup of field automorphisms is a normal subgroup of $Out(S)$, implies $G\simeq (S\times X).F$, where $F$ is a some cyclic $p$-group. In particular, we get that $X.F$ is a $p$-group. We can assume that $F=\langle xX\rangle$. As noted above $C_X(x)<X$; a contradiction with fact that $|a^G|_p=1$ for each $a\in X$.
Let $x\in G$ be a $p$-element such that $Ind_G(x)_p>1$. We have $x\in SC_G(S)$. Since $S$ and $C_G(S)$ are normal subgroups of $G$ with trivial intersection, we obtain $x$ has unique expression in form $x=x_Sx_C$ where $x_S\in S, x_C\in C_G(S)$. Moreover from Lemma \ref{Dol} it follows that $C_G(x)=C_G(x_S)\cap C_G(x_C)$.
From Lemma \ref{Simple} it follows that $S$ contains $h$ such that a $p$-element $y\in C_G(h)$ iff $y\in C_G(S)$. Therefore, for each $A\in Syl_p(C_G(h))$ there exists $B\in Syl_p(C_G(x_S))$ such that $A<B$. Since $Ind_S(x_S)_p<Ind_S(h)_p$ and $x_C\in C_G(h)$, we obtain $Ind_G(x)_p<In_G(hx_C)$; a contradiction.
\end{proof}
Let $O=O_p(G)$. Lemma \ref{cokp} implies that $O$ includes the socle of $G$. Therefore, $C_{G}(O)=Z(O)$, and for each $h\in G\setminus O$ we have $Ind_{G}(h)_p=|G||_p$.

\begin{lemma}\label{pchast}
$|G|_p>|O|$
\end{lemma}

\begin{proof}
Assume that $|G|_p=|O|$. Let $x\in O$ such that $Ind_O(x)>1$ and $h\in G$ be a $p'$-element. We have $|C_O(x)|>|Z(O)|$. Therefore, $|C_O(h)|>|Z(O)|$, consequently $C_O(h)$ contains an element $y$ such that $Ind_O(y)>1$. Hence $C_O(y)=C_O(h)$. From Lemma \ref{Thom} it follows that $C_O(h)=O$; a contradiction.
\end{proof}

 From Lemma \ref{pchast} if follows that $|G|_p>|O|$. Let $h\in G$ be a $p'$-element and $x\in C_G(h)\setminus O$ be a $p$-element. Using Lemma \ref{Thom} we can show that $Ind_G(x)_p=1$, so $x\in C_G(O) = Z(O)$, which is a contradiction. In particular $p$ is a connected component of $GK(G/O)$.

\begin{lemma}\label{ab}
The group $O$ is abelian.
\end{lemma}
\begin{proof}
Assume that there exists $x\in O\setminus Z(O)$. Put $h\in G$ is a $p'$-element. We have $|C_G(h)|_p=|C_G(x)|_p$, in particular $C_G(h)$ contains a $p$-element $y$ such that $Ind_G(y)_p>1$. From Lemma \ref{Thom} it follows that $O<C_G(h)$; a contradiction.
\end{proof}
\begin{lemma}\label{zo}
$Ind_G(x)_p=1$ for each $x\in O$.
\end{lemma}
\begin{proof}
Assume that there exists $x\in O$ such that $Ind_G(x)_p>1$. Lemma \ref{ab} implies that $O$ is abelian. Therefore, $C_G(x)\geq O$. Put $h\in G$ is a $p'$-element. We know that $p$ is a connected component of $GK(G/O)$, hence for each Sylow $p$-subgroup $P$ of $C_G(h)$ we have $P\leq O$. Since $C_G(O)=Z(O)$, we see that $P\neq O$. Consequently $|G|_p/|C_O(x)|\geq |G|_p/|C_G(x)|_p=|G||_p$. Therefore, $|P|\geq|O|$; a contradiction.
\end{proof}

\begin{lemma}\label{aga1}
The group $G$ is non solvable.
\end{lemma}
\begin{proof}
Assume that $G$ is solvable. From Lemma \ref{GK} it follows that $G/O$ is a Frobenius or $2$-Frobenius group. Since kernel of $G/O$ is a $p'$-group and Lemma \ref{norm}, we obtain $G/O$ is a Frobenius group with $p'$-kernel $\overline{K}$ and complement $\overline{F}$, else $G$ is not minimal. Put $K<G$ be a minimal subgroup such that $KO/O~=~\overline{K}$. From Frattini argument we have $N(K)O/O\simeq \overline{K}$. Let $F\leq N_G(K)$ be a minimal subgroup such that $FO/O=\overline{F}$. Since $G/O$ is a Frobenius group with the complement $\overline{F}$, we see that $N_G(F)<OF$; in particular $\pi(N_G(F))=\{p\}$.

Let $H<K$ be maximal with respect to inclusion subgroup of $K$ such that $C_O(H)>Z(G)$. We show that $H$ is not trivial. For each $h\in K$ and $y\in F\setminus O$ we have $|h^G|_p=|y^G|_p$ and $|C_G(y)|_p>|Z(G)|_p$. Therefore, $C_O(h)>Z(G)$, in particular $H>1$.

 Let $x\in C_O(H)\setminus Z(G)$. Lemma \ref{zo} implies that $x\in Z(P)$ for some Sylow $p$-subgroup $P$ of $G$. Hence $C_G(x)$ includes a subgroup $V$ which is conjugated with $F$. We have $H<C_G(x)$. Since $H$ is a maximal $p'$-subgroup with a non trivial centralizer in $O$, we see that $H$ is a Hall $p'$-subgroup of $C_G(x)$. In particular $C_G(x)= O.\overline{H}.\overline{F}$, where $\overline{H}=HO/O$. Let $y\in C_O(H)\setminus Z(G)$. We can show that $C_G(y)=O.\overline{H}\leftthreetimes \overline{R}$. Since $\overline{R}<N_{G/O}(\overline{H})$ and $\overline{F}<N_{G/O}(\overline{H})$, we have $\overline{H}\overline{R}\overline{F}$ is a Frobenius group with the kernel $\overline{H}$ and the complement $\overline{F}$. Therefore, $\overline{H}\overline{R}\overline{F}=\overline{H}\overline{R}=\overline{H}\overline{F}$ and $C_G(y)=C_G(x)$. Thus, $C_G(x)=C_G(C_O(H))$.

Let $N=N_K(H)$. Since $K$ is nilpotent, we get that $H=K$ or $N>H$. We have $N<N_G(C_G(H))$. Therefore, $N<N_G(C_G(C_O(H)))$. If $N>H$, then $N\not\leq C_G(C_O(H))$. Thus, $N_G(F)$ includes a subgroup $X$ such that $(X/H)\simeq N/H$; a contradiction with fact that $N_G(F)$ is a $p$-group. It follows that $H=K$. We have $F<C_G(C_O(H))$ and $O$ is abelian. Therefore, $C_O(H)\leq Z(G)$; a contradiction.
\end{proof}

From Lemmas \ref{aga1} and \ref{GK} it follows that $G/O\simeq F.S$, where $F$ is a nilpotent $\pi_1(G/O)$-group, and $S$ is a nonabelian simple group. Let $g\in G$ be a $p'$-element. Put $H_g\leq C_G(g)$ a subgroup generated by all $(\pi(|g|)\cup\{p\})'$-elements. We have $H_g\leq C_G(C_O(g))$. Therefore, $H_g$ acts regularly on $O/C_O(g)$. Hence Sylow subgroups of $H_g$ are cyclic or quaternion groups.

\begin{lemma}\label{aga2}
$F=1$.
\end{lemma}
\begin{proof}

Assume that $|\pi(F)|>1$. Let $h\in F$ be a $t$-element for some $t\in\pi(F)$. We have $C_F(h)$ includes some Hall $t'$-subgroup $T$ of $F$. Since Sylow subgroups of $H_h$ are cyclic or quaternion, we get $S$ acts trivially on $T$. Therefore, $F$ is a $t$-group for some prime $t\neq p$.

Assume that there exists $h\in Z(F)$ such that $|\pi(C_{F.S}(h))|>1$. Let $g\in C_{F.S}(h)$ be a $t'$-element. Assume that $Z(F)<C_{F.S}(g)$. Since $Z(F)$ is a normal subgroup of $F.S$, we have $\langle g^G\rangle<C_{F.S}(Z(F))$. The group $S$ is a simple group, therefore, $\langle g^G\rangle F/F=S$. In particular $F.S$ contains an element of order $pt$; a contradiction. Hence $[Z(F),g]>1$. We have $H_h$ acts regularly on the $O/C_O(h)$ and $\langle g^{C_{F.S}(h)}\rangle\leq H_h$. This assertion contradicts with Lemma \ref{l:hallhigman}. Therefore, $\pi(C_{F.S}(h))=\{t\}$ for each $h\in Z(F)$.

Assume that there exists $a\in F.S$ a $t$-element such that $|\pi(C_{F.S}(a))|>1$. Since $C_{F.S}(a)\cap Z(F)>1$, we get a contradiction with Lemma \ref{l:hallhigman}.

Therefore, $t$ is a connected component of $GK(F.S)$. Since $t\in \pi_1(GK(F.S))$, we get $t=2$. From a description of the prime graph of finite simple groups \cite{VasilevVdovin} it follows that $S\simeq Alt_5$. From Brauer $2$-character tables \cite{Conway} it follows that $F.S$ contains an element $x$ such that $|x|\in\{6,10\}$; a contradiction.
\end{proof}

From Lemmas \ref{aga1} and \ref{aga2} it follows that $G/O$ is a simple group and $p$ is a connected component of $GK(G/O)$.

\begin{lemma}\label{Fin}
$G\simeq O$
\end{lemma}
\begin{proof}
Assume that $S$ is an alternating simple group of degree $n$. An alternating group has disconnected prime graph iff one of the numbers $n,n-1,n-2$ is prime, and this number is a connected component of the prime graph. In particular, if $n>6$, then $3\in \pi_1(S)$. If $n>6$, then $S$ contains an element $g$ of order $3$ such that $C_S(g)\simeq \langle g\rangle\times Alt_{n-3}$. Therefore, in this case $H_g$ includes a Frobenius group; a contradiction with assertion that $H_g$ acts regular on $O/C_O(g)$.
Let $n\in\{5,6\}$. Therefore, $S$ contains an element $g$ such that for each $p$-element $h\in S$ we have $\langle h,g\rangle=S$. We have $C_O(g)>Z(G)$. Put $x\in C_O(g)\setminus Z(G)$. It follows from Lemma \ref{ab} that $|x^G|_p=1$. Hence, $C_G(x)$ includes a Sylow $p$-subgroup $P$ of $G$. In particular $\langle g,PO/O\rangle=S$. That signifies that $C_G(x)=G$; a contradiction.

Assume that $S$ is a group of Lie type. If Lie rank of $S$ is more then $2$, then $S$ contains an element $g$ such that $H_g$ includes a Frobenius group. Therefore, we can assume that Lie rank of $S$ is $1$ or $2$. Assume that $S\simeq L_2(q)$. We have that $S$ is generated by a pair $a,b$ where $|a|=(q+1)/2, |b|=(q-1)/2$. Since $p=(q-1)/2$ or $p=(q+1)/2$, we can assume that $C_O(g)\leq Z(G)$; a contradiction. Groups $L_3(q)$ and $U_3(q)$ contain an element $g$ such that $C_S(g)$ includes $L_2(q)$. Therefore, $S$ is not isomorphic to one of a $L_3(q)$ or $U_3(q)$. Similarly, it can be shown that $S$ is not isomorphic to $B_2(q), ^2B_2(q), G_2(q), ^2G_2(q)$ and sporadic groups.
\end{proof}

Lemma \ref{Fin} completes proof of the theorem.

\section{Proof of Corollaries}
Proof of Corollary 1.
\begin{proof}
If $G$ is a $p$-group, then the corollary is satisfied. Let $h\in G$ be a $p'$-element. We have $|C_G(h)|_p>|Z(P)|$, where $P\in Syl_p(G)$. Let $H\in Syl_p(C_G(h))$. It follows from Lemma \ref{pshtrih} that $G/O_{p'}(G)\in R^*(p)$. From Theorem 1 we get that $G/O_{p'}(G)$ is a $p$-group. Since $|H|>|Z(P)|$, we get that $H$ contains $x$ such that $Ind_{G/O_{p'}(G)}(xO_{p'}(G))=p^e$, for some $e>0$. Therefore, $Ind_{G/O_{p'}(G)}(xO_{p'}(G))_p=Ind_G(x)_p=p^e$. Hence $C_G(h)$ includes some Sylow $p$-subgroup of $C_G(x)$. From Lemma \ref{hh} it follows that $Z(P)\leq Z(G)$.
\end{proof}

Proof of Corollary 2.

\begin{proof}
If $G\in R^*(p)$, then from Corollary 1 it follows that $Z(G)>1$; a contradiction. Therefore, $Ind_G(x)_p=1$ for each $p$-element $x$ of $G$. From Lemma \ref{Navar} it follows that $O^{p'}(G/O_{p'}(G))= S_1\times\cdot\cdot\cdot \times S_r\times H$. From Lemmas \ref{norm} and \ref{pshtrih} it follows that $O^{p'}(G/O_{p'}(G))\in R(p)$. Since $S_1$ is a subnormal subgroup of $G/O_{p'}(G)$, we get $|x^{G/O_{p'}(G)}|$ is a multiple of $|x^{S_1}|$ for each $x\in S_1$. Assume that $p=3$ and $S_1\simeq Ru$. We have $S_1$ contains $x$ of order $12$ such that $|C_{S_1}(x)|=24$. Therefore, $1<|x^{S_1}|_3<|G||_3$; a contradiction. If $S_1\simeq J_4$, then $S_1$ contains $x$ of order $21$ such that $|C_{S_1}(x)|=42$; a contradiction with definition of $R^*(p)$-groups. Assume that $S_1\simeq ^2F_4(q)'$. According to the results of \cite{Shinoda}, there is just one conjugacy class of
elements of order $3$ in $S_1$. Therefore, $N_{S_1}(\langle x\rangle)\simeq C_{S_1}(x):2=3.U_3(q):2$ where $x\in S_1$ is an element of order $3$ \cite{Malle}. Thus, $S_1$ contains an element $y$ of order $6$ such that $1<|y^{S_1}|_3<|G||_3$; a contradiction.
Assume that $p=5$ and $S_1\simeq Th$. In this case $S_1$ contains an element $x$ of order $8$ such that $|C_{S_1}(x)|=96$; a contradiction.
The assertion of Corollary 2 follows from Lemma \ref{Navar}.
\end{proof}

\EditInfo{%
    October 03, 2019}{%
    March 04, 2020}{%
    Ivan Kaygorodov}

\end{paper}